# Math goes to Hollywood
## Stereotypen in Filmen und Serien dekodieren

András Bátkai & Ingrid Gessner

*Stereotype Vorstellungen von Mathematik und Mathematiker\*innen beeinflussen das Interesse von Jugendlichen an MINT-Fächern. Daher plädiert dieser Beitrag dafür, populäre Filme und erfolgreiche Serien nicht nur im Fremdsprachenunterricht, sondern auch im Mathematikunterricht einzusetzen. Durch die Analyse audiovisueller Medien im Unterricht können verzerrte Darstellungen über Mathematiker\*innen erkannt, die Gründe dafür benannt und alternative Sichtweisen entwickelt werden. Dies trägt dazu bei, das Bild der Mathematik zu entmystifizieren und den Spaß an Mathe und Englisch zu fördern. Die Einbindung sozialwissenschaftlicher Themen steigert zudem die Motivation der Schüler\*innen.*

**Schlagwörter:** Filme im Mathematikunterricht, Filme im Fremdsprachenunterricht, Fächerübergreifender Unterricht, Vorurteile abbauen, Ideologiekritik, Reflexive Geschlechterpädagogik

## Einleitung

Populäre Filme, Serien und Videoclips im Fremdsprachenunterricht zu nutzen ist lange erprobt und hat nicht zuletzt durch die digitale Verfügbarkeit zu einer verstärkten Nutzung im Englischunterricht geführt. Audiovisuelle Medien eignen sich nicht nur zur praktischen Spracharbeit, zur Thematisierung kultureller Bezüge, sondern auch zur medienkritischen Auseinandersetzung mit den vermittelten Inhalten und werden im Common European Framework (CEFR, 2018) explizit für den Einsatz im Fremdsprachenunterricht empfohlen. Dass die Nutzung von Blockbustern und Serienhits nicht auf den Fremdsprachenunterricht beschränkt bleiben sollte, sondern auch im Mathematikunterricht sinnvoll und zielführend ist, dafür möchte dieser Beitrag eine Lanze brechen.

In den letzten 25 Jahren hat Hollywood eine Vielzahl von Filmen mit Mathematiker\*innen in Haupt- oder Nebenrollen produziert. Was haben diese Filme gemeinsam? Eine Kreidetafel oder ein Whiteboard voller Zahlen und Formeln, das gleiche auf Servietten und losen Blättern und die Person, von der diese Notizen stammen ist unterschätzt, in sich gekehrt, nerdig bis klinisch autistisch. Der oder die Protagonist\*in besticht gleichzeitig durch Intelligenz, Einfallsreichtum, Neugierde und Engagement. Mit wenigen Ausnahmen ergeben diese Zutaten den klassischen Blockbuster oder Serienhit über Mathematiker\*innen. Die Hauptrollen sind meist männlich, aber auch immer wieder weiblich besetzt.

In diesem Beitrag interessieren uns die Offenlegung von Klischees in Filmen über Mathematiker\*innen und die interdisziplinären Überlappungen von Mathematik, Sprach-, Kultur- und Literaturwissenschaften in populären Filmen und Serien über Mathematiker\*innen. Wir möchten darüber hinaus dazu anregen, Stereotypen über Mathematiker\*innen im Unterricht zu dekodieren, um dazu beizutragen, das Bild der Mathematik zu entmystifizieren und damit nutzbar zu machen, um „Spaß an Mathe" und „Spaß an Englisch" zu fördern. Hierzu kann auf eine Vielzahl bekannter Ansätze und Methoden aus der Filmdidaktik zurückgegriffen werden (Henseler, Möller & Surkamp, 2011; Thaler, 2012; Viebrock, 2016). Die Lernenden erarbeiten die audiovisuellen Texte lernerzentriert, untersuchen die unterschiedlichen Bedeutungsebenen und formulieren verschiedene Interpretationen.

Im Englischunterricht könnte darüber hinaus die filmische Repräsentation und ggf. Kritik an traditionellen amerikanischer Werten wie Individualismus, das Erreichen und Überwinden der (Western) Frontier und der amerikanische Traum (American Dream) im Mittelpunkt stehen. Einen Mehrwert haben diese Erkenntnisse aber auch im Mathematikunterricht. Denn durch diesen multidisziplinären Ansatz können die Schüler\*innen ein tieferes Verständnis dafür entwickeln, wie in der Populärkultur (amerikanische) Werte, Ideologien und Stereotypen unser Verständnis von Mathematik als Fach prägen und verstärken. Darüber hinaus können die Schüler\*innen durch die Erforschung der mathematischen Konzepte und Techniken, die in diesen Filmen verwendet werden, Einblicke in die Rolle der Mathematik in unserer Gesellschaft erhalten. Die Schüler*innen entwickeln ihre Fähigkeiten zum kritischen Denken weiter und können die Bedeutung der Mathematik in unterschiedlichen Bereichen und Kontexten erkennen und bewerten.

Dies deckt sich mit dem österreichischen Lehrplan (BMBWF, 2023) für die Unterstufe, der zum übergreifenden Thema „Interkulturelle Bildung" vorsieht, dass Schüler*innen befähigt werden, „Stereotype, (Fremd-)Zuschreibungen und Klischees [zu] identifizieren sowie ausgrenzende, rassistische, sexistische Aussagen und Handlungsweisen erkennen, hinterfragen und dagegen auftreten." Reflexive Geschlechterpädagogik und Fragen der Gleichstellung sind Themen sowohl im Englisch- als auch im Mathematikunterricht. Auch Heymann zählt zu den zentralen Aufgaben der auf Allgemeinbildung abzielenden Schule (Heymann, 1996) die Verbindung des zu unterrichtenden Material und dessen Anwendungen mit anderen Interessens- und Lebensbereichen.

## Wahrnehmungen

Ernest (1998) beschreibt das in der Öffentlichkeit verbreitete Bild der Mathematik mit den Adjektiven schwierig, kalt, abstrakt, theoretisch und ultra-rational. Mathematik wird darüber hinaus als sehr wichtig und weitgehend männlich wahrgenommen. Wenn man aber nachfragt, warum wäre Mathematik wichtig, bekommt man meistens unsichere Antworten. In Film und Fernsehen, Internet, Videospielen und Büchern (Massenmedien) wird die Mathematik zumeist als ein schwieriges Fachgebiet dargestellt, das nur einigen wenigen zugänglich ist, die zudem in einer genderspezifischen und negativen Weise dargestellt werden (Hall & Suurtamm, 2020). Mathematiker*innen auf der Leinwand und in Serien sind introvertiert, sozial beeinträchtigt und ihre mathematische Genialität wird oft mit Verrücktheit oder psychischer Labilität gleichgestellt (Latterell & Wilson, 2004; Wilson & Latterell, 2001) während die Darstellung dessen, was Mathematik bedeutet auf die Berechnung komplizierter Formeln beschränkt wird.

Wittmann (2003) beoabachtet: „Das dem traditionellen Mathematikunterricht zu Grunde liegende Mathematikbild erschwert oder versperrt vielen Lernenden den Zugang zur Mathematik." Er stellt aber weiterhin fest, dass mathematische Inhalte oft als streng geordnete Systeme von klaren Begriffen, Regeln und Verfahren wahrgenommen werden, die speziell für bestimmte Aufgaben entwickelt wurden. Der Lernprozess erfolgt schrittweise und portionsweise, um eine möglichst fehlerfreie Reproduktion zu erreichen. Die Unterscheidung zwischen dem, was „wahr" und „falsch" ist, ist scharf definiert, wodurch die Angst vor inhaltlichen oder formalen Fehlern groß ist. Viele Lernende glauben, dass mathematische Fähigkeiten erforderlich sind, um Mathematik verstehen zu können. Aus diesem Grund werden eigene Überlegungen, die zunächst fehlerhaft sein können, oft zugunsten der Reproduktion vorgegebener Musterlösungen vernachlässigt. Aufgrund eines Mangels an Verständnis klammern sich viele Lernende an die äußere Form der Darstellung, hinter der die eigentlichen Inhalte oft verborgen bleiben.

Chronaki und Kollosche (2019) haben Wittmanns Beobachtungen empirisch untermauert. Schüler*innen, die sich von Mathematik abwenden, tun dies aus den von Wittmann beschriebenen Gründen, nämlich aus Angst, Fehler zu machen und der Überzeugung, nicht (genügend) begabt zu sein. Für viele kommt dazu, dass man oft allein arbeiten muss und dass Mathematik mit Rechnen gleichgestellt wird.

Auch wenn es nicht nur die Mathematik betrifft, muss uns zu denken geben, dass stereotype Vorstellungen von MINT-Fachkräften sehr wahrscheinlich das Interesse von Mädchen und jungen Frauen beeinflussen, sich in diesem Berufsfeld qualifizieren zu wollen (Buelvas-Baldiris & Rubira-García, 2023; Steinke & Paniagua Tavarez, 2018). Es ist umso alarmierender, dass die entsprechenden Vorurteile aus Wahrnehmungen herrühren, die aus Film- und Seriendarstellungen von Wissenschaftler*innen entnommen werden oder durch diese verstärkt werden. Nicht nur Schüler*innen, sondern auch angehende Lehrkräfte haben oft sowohl falsche Vorstel-lungen von der Mathematik als Wissenschaft (Schreck, Groß Ophoff & Rott, 2023) als auch von Mathematiker*innen (Woltron, 2020).

In diesem Zusammenhang sei angemerkt, dass eine groß angelegte transnationale Studie der EU (European Commission, 2022) festgestellt hat, dass Bildungssysteme, die sozialwissenschaftliche Themen in ihren Lehrplänen verankert haben, einen höheren Anteil von 15-jährigen Schüler*innen aufweisen, die über grundlegende wissenschaftliche Kompetenzen verfügen. Diese Erkenntnis ist für die Didaktik der Mathematik von Interesse, da sie zeigt, dass die Behandlung medienwissenschaftlicher und sozialkritischer Themen dazu beitragen kann, die Mathematik zu entmystifizieren und die Motivation zu steigern.

In ihrer Studie zu Bilderbüchern mit mathematischen Themen haben Fellus et al. (2022) festgestellt, dass selbst Geschichten, die bewusst Stereotype über Fähigkeiten, Geschlecht und ethnische Zugehörigkeit ablehnen, andere Vorstellungen darüber verstärken können, was Mathematik ausmacht und was es bedeutet, ein „Mathe-Mensch" zu sein. In zahlreichen Bilderbüchern wird gezeigt, dass mathematisch begabte Menschen eine fast übernatürliche Begabung für schnelles und fehlerfreies Rechnen besitzen und/oder Schwierigkeiten haben, soziale Beziehungen aufzubauen.

Ein wesentliches Ziel der kritischen Analyse populärer Filme und Serien im Unterricht muss es deshalb sein, Bedeutungsebenen zu untersuchen, Vorurteile zu erkennen und einzuordnen, um dadurch zu einem realistischeren Bild von Mathematik und Mathematiker*innen zu gelangen. Wenn im folgenden Abschnitt Filme und Serien vorgestellt werden, dient dies zum einen dazu das weite Feld der Produktionen aufzuzeigen und zum anderen dazu das Augenmerk auf die problematische Darstellung von Mathematik und Geschlecht zu lenken. Der implizierte nächste Schritt, nämlich die Medienproduktionen in ihrer Gesamtheit oder selektiv in Unterrichtssequenzen didaktisch aufzubereiten, diese durchzuführen und zu evaluieren, kann in diesem Beitrag nur stellenweise mit einigen Hinweisen geleistet werden. Für die Aufbereitung sei nochmal auf zahlreiche Publikationen in der Englisch- und Filmdidaktik verwiesen (Bradley, 2016; Henseler et al., 2011; Lütge, 2012; Thaler, 2014; Viebrock, 2016).

## Math goes to Hollywood

### *The Imitation Game – Ein streng geheimes Leben*

Einer der bekanntesten Mathematiker*innenfilme des letzten Jahrzehnts ist *The Imitation Game – Ein streng geheimes Leben* (2014) mit Benedict Cumberbatch als Alan Turing, dem britischen Mathematiker, der während des Zweiten Weltkriegs eine entscheidende Rolle beim Knacken des Enigma-Codes spielte. Der Film porträtiert Turing als brillanten Mathematiker, dessen Fähigkeiten für die militärischen Ziele und den Sieg der Alliierten über Nazi-Deutschland unerlässlich waren. Gleichzeitig bedient *The Imitation Game* allerdings auch das Klischee des nicht zu persönlichen Beziehungen fähigen, verkannten Genies, dem Zahlen näher sind als Menschen. Im Film scheitert Turing an seiner Identität als Ausgegrenzter, als beziehungsunfähiger Wissenschaftler und schwuler Mann in einer Gesellschaft, die allem Andersartigen gegenüber feindlich eingestellt ist.

An der Seite Cumberbatchs als Turing spielt Keira Knightley die Mathematikerin und Kryptoanalystin Joan Clarke. Wer glaubt, dass die Figur der Mathematikerin einen Schritt in Richtung Gleichberechtigung bedeutet oder Joan Clarke als Identifikationsfigur dienen könnte, wird enttäuscht. *The Imitation Game* besteht den Bechdel-Test, der oft verwendet wird, um die Präsenz von Frauen in Filmen und anderen Medien zu bewerten, nicht. Um den Test zu bestehen, muss ein Film drei Kriterien erfüllen: 1) der Film muss zwei namentlich genannte Frauen haben (2), die sich miteinander unterhalten (3) und zwar über ein Thema, das nicht Männer betrifft. Die Regeln des Bechdel-Tests können von der Lehrperson erläutert werden und von den Schüler*innen diskursiv auf den Film bzw. Sequenzen wie die Verlobungsfeier von Turing und Clarke angewendet werden.

Clarke hat ihren Abschluss in Cambridge gemacht, hatte aber als Frau keine Chance, ihre Studien dort fortzuführen. Sie taucht im Film als Anwärterin für das Team um Turing auf, nachdem dieser ein Rätsel in einer Zeitungsannonce veröffentlicht hat, um neue Mitarbeiter*innen zu finden. Als sie am angegebenen Ort eintrifft, wird sie sofort von einem Mann aufgehalten, der glaubt, dass sie sich um eine Stelle als Sekretärin bewirbt, da sich vor ihrer Ankunft nur Männer gemeldet hatten. Auch wenn sie anschließend wegen ihrer hervorragenden Arbeit eingestellt wird, um in Bletchley Park and der Seite von Turing zu arbeiten, beträgt Knightleys Leinwandzeit nur einen Bruchteil derjenigen von Cumberbatch. Die Figur der Joan Clarke dient ausschließlich Turings Charakterentwicklung und nicht dem Vorantreiben der Handlung oder dem Lösen mathematischer Probleme.

### *Good Will Hunting*
Auch in *Good Will Hunting* (1997) wird der Mathematiker als gestörtes Genie dargestellt, das in diesem Fall mit seiner kriminellen Vergangenheit und verdrängten Missbrauchserfahrung hadert. Die Hauptfigur Will Hunting (gespielt von Matt Damon) wächst in einfachen Verhältnissen in Boston auf, wo er ironischerweise als Putzmann und nicht als Wissenschaftler an der renommierten Universität MIT (Massachusetts Institute of Technology) arbeitet. Der Film zeichnet nach, wie Hunting sich mit der Unterstützung seines Therapeuten, gespielt von Robin Williams, seiner Vergangenheit stellt. Er lernt, seine Fähigkeiten zu akzeptieren und zu nutzen und ergreift am Ende sogar die Initiative, um seine Liebesbeziehung zu retten. *Good Will Hunting* folgt der Formel „Vom-Tellerwäscher-zum-Millionär" („from rags to riches") bzw. der Aschenputtel-Geschichte und hat als modernes Mathematikmärchen nur wenig mit der Realität zu tun. Deshalb taugt Will Hunting nur bedingt als mathematische Identifikationsfigur für Menschen, die in ähnlichen sozioökonomischen Verhältnissen wie er aufgewachsen sind. Die Mathematik ist, wie der Prinz im Märchen, nur wenigen Auserwählten vorbehalten. Die Absurdität dieser Aussage ist leicht von den Schüler*innen zu durchschauen und kann im Unterricht diskutiert und entkräftet werden.

### *Proof – Der Beweis: Liebe zwischen Genie und Wahnsinn* and *An Invisible Sign*
Bei *Proof – Der Beweis: Liebe zwischen Genie und Wahnsinn* (2005) zeigt bereits der deutsche Untertitel auf, worum es neben der Mathematik vor allem geht: Genialität gekoppelt mit Wahnsinn, verkompliziert durch eine Liebesgeschichte. Der Film zeichnet die Vater-Tochter-Beziehung zweier mathematischer Genies nach. Während der Vater im Alter von 27 Jahren aufgrund einer schweren Krankheit ins Krankenhaus eingeliefert wird, beschließt die Tochter Catherine Jahre später, aber ebenfalls mit 27 Jahren, ihre Universitätskarriere abzubrechen, um bei ihrem Vater zu leben und zu sich selbst zu finden. Aus Sorge, dass sie sowohl die Intelligenz als auch die Krankheit ihres Vaters erben könnte, beginnt Catherine Halluzinationen zu haben, die nach dem Tod des Vaters zunehmen. In *Proof* wird also erneut das Klischee von mathematischer Genialität und Geisteskrankheit ausgespielt. Einem ganz ähnlichen Muster folgt der Film *An Invisible Sign* (2010). Auch hier ist es eine junge Frau, Mona Gray (gespielt von Jessica Alba), die vollkommen zurückgezogen nur in der Welt ihrer Zahlen und der Mathematik lebt, während ihr Vater, ebenfalls Mathematiker, psychisch erkrankt ist.

Dass die Gleichung von Genie und Wahnsinn fehlerhaft ist, kann für beide Filme leicht herausgearbeitet werden, indem die Lehrperson zu bedenken gibt, dass nicht jeder kreative oder mathematisch begabte Mensch geisteskrank und umgekehrt nicht jeder Geisteskranke kreativ oder mathematisch begabt ist. Anschließend können die folgenden Fragen gestellt werden: Inwiefern werden Stereotype oder Klischees über Genies und Mathematiker*innen in den Filmen aufgegriffen? Wie können sie widerlegt werden?

### *Agora – Die Säulen des Himmels*
*Agora – Die Säulen des Himmels* (2009) ist die Hollywoodversion des Lebens der spätantiken Mathematikerin Hypatia, von der heute oft nicht mehr als deren spektakuläre Ermordung durch fanatische Christen bekannt ist. Die Astronomin, Mathematikerin und Philosophin Hypatia wurde um das Jahr 370 n.Chr. in Alexandria geboren, wo sie am Museion, der führenden For-schungseinrichtung jener Zeit, unterrichtete. Sie interpretierte die Lehren Platons und Aristoteles' und die Werke des Astronomen Ptolemäus, der drei Jahrhunderte vor ihr lebte. Auch wenn von ihren Schriften nichts überliefert ist, wissen wir aus den Briefen und Schriften des Synesios von Kyrene, dass sie zur Verbesserung des Astrolabiums, einem Instrument zur Bestimmung der Position von Himmelskörpern, beitrug und die erste Sternenkarte in Form einer Planisphäre weiterentwickelte. Sokrates Scholastikos schreibt in seiner *Kirchengeschichte*: „Es gab in Alexandria eine Frau mit Namen Hypatia, Tochter des Philosophen Theon, die in Literatur und Wissenschaft so erfolgreich war, dass sie alle Philosophen ihrer Zeit übertraf."

In *Agora* spielen Hypatias philosophische und wissenschaftliche Arbeiten eine Rolle, aber im Vordergrund stehen Hypatias Schönheit und die Männer, die in sie verliebt waren; allen voran die historisch nicht belegte Figur des Sklaven Davus, der Hypatia im Film erstickt, um sie von einem noch grausameren Tod zu bewahren. Dies ändert nichts an der Tatsache, dass die nicht dem Idealbild der antiken Frau entsprechende Hypatia sterben muss – weder ein Erfolg für die Frauen, noch für die Mathematik, aber bis heute ein beliebtes literarisches Ende für Frauen, die sich nicht in das von der Gesellschaft vorgegebene Frauenbild einfügen. In Unterricht sollte die Brisanz dieser Erkenntnis für Diskussionsstoff sorgen.

### *Hidden Figures – Unerkannte Heldinnen*
*Hidden Figures – Unerkannte Heldinnen* (2016) ist da ein positiveres Beispiel für einen erfolgreichen Film mit Mathematikerinnen in den Hauptrollen. Der Film, der sich erzählerische Freiheiten nimmt, dokumentiert die auf wahren Gegebenheiten beruhende Geschichte einer Gruppe afroamerikanischer Mathematikerinnen, die während des Wettlaufs um die erste Rakete im All in den 1950er und 1960er Jahren für die NASA arbeiteten und Berechnungen von Hand durchführten. Für die Behandlung dieses Films im Unterricht ist es wichtig, dass die Lehrkraft den Schüler*innen das notwendige Hintergrundwissen zur Verfügung stellt. So steht 1962 der Kalte Krieg auf seinem Höhepunkt. Im gleichen Jahr hatte Präsident John F. Kennedy den Weltraum, bzw. den Mond als die New Frontier definiert, die es bis zum Ende des Jahrzehnts zu erreichen gilt. Er bezog sich dabei auf die Western Frontier, die bereits Ende des 19. Jahrhunderts als erreicht bzw. geschlossen („closed") erklärt worden war. Der Historiker Frederick Jackson Turner nahm dies 1893 zum Anlass, die sogenannte Frontier-Hypothese zu formulieren, die postuliert, dass der genuin amerikanische Charakter an der Grenze geformt wurde und wird. Dass dies nicht der Fall ist, ist seither zahlreich bewiesen worden, der Mythos hat sich allerdings gehalten und findet sich nicht nur in Filmen wie *Hidden Figures*.

Im Film geht es darum, dass vom Präsidenten proklamierte Ziel der New Frontier zu erreichen. Ein wenig spielt aber auch die Ideologie des amerikanischen Traums eine Rolle, wenn Katherine Johnson (gespielt von Taraji P. Henson) entscheidende analytische Daten für das Projekt liefert, das den US-Amerikaner John Glenn in die Erdumlaufbahn bringt, während Mary Jackson (Janelle Monáe) an der Konstruktion von Glenns Rakete mitwirkt. Dorothy Vaughan (Octavia Spencer) bringt sich selbst die Programmiersprache FORTRAN bei und wird zur Leiterin der IBM-Maschinengruppe am NASA-Standort Langley. Dabei erfahren alle drei Diskrimination, und zwar nicht nur, weil sie Frauen sind, sondern weil sie schwarze Frauen sind. Umgekehrt beindrucken sie Zuschauer*innen und Vorgesetzte durch ihre Intelligenz und Entschlossenheit angesichts von Sexismus und Rassismus.

### *Inception* und *Numb3rs*
In *Inception* (2010) nutzt ein Team von Spezialist*innen um das charismatische Ehepaar Dominick bzw. „Dom" und Mal Cobb (gespielt von Leonardo di Caprio und Marion Cotillard) ihre Fähigkeiten in verschiedenen Bereichen, darunter neben Pharmakologie auch Mathematik, um einen anspruchsvollen Raub zu begehen. Der Mathematiker in der Gruppe, Arthur (gespielt von Joseph Gordon-Levitt), ist dafür verantwortlich, die physikalischen Gesetze in der Traumwelt zu manipulieren, um dem Team zu helfen, dessen Ziele zu erreichen. In einer der spektakulärsten Action-Szenen des Films gelingt es Arthur die Schwerkraft in einem Hotelflur so zu verändern, dass die Teammitglieder in der Luft schweben und gegen Wände laufen können, um so eine feindliche Gruppe von Traum-Manipulatoren abzuwehren und gleichzeitig die Träumenden, die sie infiltriert haben, vor Entdeckung zu schützen. Die Mathematik wird zum Werkzeug sinisterer Absichten missbraucht. Dennoch ist man als Zuschauer*in auf der Seite der kriminellen Wissenschaftler*innen, weil das eigentliche Ziel – ganz im Sinne der Trope des letzten Verbrechens – es sein soll, dem Hauptdarsteller durch die Tilgung dessen Strafregisters die Rückkehr zu seinen Kindern zu ermöglichen. Ob dies tatsächlich gelingt oder auch nur eine weitere Traumsequenz ist, bleibt am Schluss offen bzw. der Interpretation des Publikums überlassen.

Während mathematische Fähigkeiten in *Inception* für Verbrechen missbraucht werden, befindet sich die Mathematik in der Fernsehserie *Numb3rs* (2005-2010) auf Seiten von Recht und Gesetz. In der Serie hilft Charlie Eppes (gespielt von David Krumholtz), ein Mathe-Genie und Professor an der fiktiven CalSci-Universität in Kalifornien, seinem Bruder Don, einem FBI-Agenten, bei der Aufklärung von Kriminalfällen. Es sind dabei vor allem Charlies ungewöhnliche Methoden und kreative Lösungsansätze, die zur Ermittlung von Tätern führen. Die Mathematik ist dabei durch die Präsenz der Mathematikerin Dr. Amita Ramanujan (gespielt von Navi Rawat) keine rein männliche Domäne. Auch sie hilft als Beraterin für das FBI bei der Lösung komplexer Fälle und sorgt dafür, dass mathematisches Wissen auch weiblich besetzt ist. Hierarchisch tritt die als ehemalige Studentin Charlies und Love Interest in die Serie eingeführte Mathematikerin hinter Charlie zurück und dient vor allem dessen Charakterentwicklung.

### *The Queen's Gambit – Das Damengambit*
Einen großen Publikumserfolg erzielte in Coronajahr 2020 die vom Streamingdienst Netflix produzierte Miniserie *The Queen's Gambit (Das Damengambit)* über ein fiktives Schachwunderkind mit außergewöhnlicher mathematischer Begabung. Das Leben von Elizabeth „Beth" Harmon (gespielt von Anya Taylor-Joy) ist geprägt vom Selbstmord der ebenfalls mathematisch hochbegabten Mutter. In einem szenischen Rückblick in Beths Erinnerung erfährt man, dass die Mutter eine Dissertation mit dem Titel „Monomial Representations and Symmetric Presentations" verfasst und an der Cornell University gelehrt hat, bevor sie sich gemeinsam mit ihrer Tochter in einen Wohnwagen in den Wäldern von Kentucky zurückgezogen hat. Der Titel der Dissertation kann als Anspielung auf das Schachspiel verstanden werden, denn eine Monomialdarstellung weist jeder Variablen ein Zahlenfeld zu, in dem die Zahlen, die nicht Null sind, wie ein maximaler Satz von Türmen auf einem Schachbrett angeordnet sind, die sich nicht gegenseitig schlagen können.

Im Waisenhaus, in dem Beth nach dem Tod der Mutter aufwächst, werden ihr und den anderen Mädchen routinemäßig Psychopharmaka verabreicht, um diese ruhig zu stellen. Zu der für die Mathematikerin im Hollywoodfilm scheinbar charakteristischen psychologischen Instabilität, kommt im Fall von *Damengambit* die im Waisenhaus erworbene Tabletten- und spätere Alkoholabhängigkeit der Hauptfigur hinzu. Trotz dieser Hindernisse schafft es Beth in die Männerwelt des Schachsports vorzudringen, sie tritt selbstbewusst auf und wird schließlich die erste amerikanische Schachweltmeisterin.

Wie bei *Hidden Figures* ist der Subtext der Serie der Kalte Krieg, wobei die Sowjetunion das Schachspiel dominiert, bis Beth 1968 von einem Vertreter des US-Außenministeriums begleitet jenseits des Eisernen Vorhangs nach Moskau reist, um am dortigen Schachturnier teilzunehmen. In der letzten Partie kommt es zum Showdown zwischen den USA und der Sowjetunion mit der amtierenden US-Meisterin und dem russischen Weltmeister und Großmeister Vasily Borgov als jeweilige Stellvertreter. Das Ideal des amerikanischen Individualismus wird durch die Serie insofern in Frage gestellt, als dass Beth, die bis zu diesem Zeitpunkt für sich alleine gekämpft hat, das Turnier nur gewinnen kann, indem sie sich gemeinsam mit ihren Freunden in den USA darauf vorbereitet. Sie kopiert damit eine Strategie, die von den überlegenen russischen Spielern seit jeher genutzt wird. Im Verlauf des Turniers wird Beths Können immer wieder mit dem der Schach-Supermacht der Sowjets verglichen. Als ihr Sieg in einer Feier vom US-Außenministerium als „die Sowjets mit ihren eigenen Waffen schlagen" interpretiert werden soll, lehnt sie ab. Stattdessen endet die letzte Folge der Serie mit einem Schachspiel, zu dem Beth Harmon von ihren russischen Fans bei einem Spaziergang in einem Moskauer Park eingeladen worden ist.

### *Arrival*

Im Film *Arrival* (2016) müssen eine Sprachwissenschaftlerin und ein Mathematiker zusammenarbeiten, um mit gliederfüßigen außerirdischen Wesen zu kommunizieren, die in zwölf urnenförmigen Raumschiffen an verschiedenen Orten der Erde gelandet sind. Sowohl die Sprachwissenschaftlerin Louise Banks (gespielt von Amy Adams) als auch der Mathematiker Ian Donnelly (gespielt von Jeremy Renner) nutzen ihr Fachwissen, um die Sprache der Außerirdischen zu entschlüsseln, deren Absichten zu verstehen und die Erde vor einem dritten Weltkrieg zu bewahren, der sich aus dem Konkurrenzdenken der Großmächte zu entwickeln droht.

Im Gegensatz zum Netflix-produzierten Film *Don't Look Up* (2021) erkennen die verantwortlichen Politiker*innen in *Arrival* glücklicherweise, dass die Auslöschung der menschlichen Existenz nur durch die konzertierte Anstrengung und Bündelung wissenschaftlicher Diskurse und Erkenntnisse abgewendet werden kann. Die anfängliche Skepsis der Mathematik gegenüber Literatur und Sprache weicht in der Aussage des Mathematikers Ian Donnelly der Anerkennung seiner Kollegin Louise Banks: „Du näherst Dich der Sprache wie ein Mathematiker."

Mit Blick auf die erzählerische Funktion stehen die beiden Hauptfiguren für die Bedeutung interdisziplinärer Ansätze bei der Problemlösung. Nur im Zusammenspiel mathematischer und sprachwissensschaftlicher Methoden können die großen Probleme gelöst werden. Ian dient als kontrastierende Figur zu Louise, die sich mehr auf die emotionalen und erfahrungsbezogenen Aspekte des Kommunikationsprozesses konzentriert, was die Bedeutung sowohl der rationalen als auch der emotionalen Intelligenz für eine erfolgreiche Kommunikation unterstreicht.

Es verwundert kaum, dass die entscheidendere Rolle in der Entschlüsselung der Sprache der Außerirdischen der Linguistin zukommt, während es der Mathematiker ist, der mit der gewissermaßen banalen Aussage, dass .083333 genau 1/12 entspricht, den entscheidenden Hinweis zur Rettung der Menschheit gibt. Dass die Sprachwissenschaftlerin Banks und der Mathematiker Donnelly sich dabei auch ineinander verlieben und am Ende des Films die gemeinsame Tochter zeugen, muss dabei wohl, als den allgemeinen Publikumserwartungen entsprechend, in Kauf genommen werden. Durch die Vorwegnahme des frühen Todes des gemeinsamen Kindes und das Scheitern der Beziehung wird ein Abdriften ins romantisch Kitschige dann auch wieder wohltuend vermieden.

## Überlegungen und Hinweise

In einer Unterrichtsreihe könnte man mit der Analyse einschlägiger Hollywood-Filme beginnen, in denen Mathematiker*innen Hauptrollen spielen, z. B. *The Imitation Game, Hidden Figures* und *Arrival*. Die Schüler*innen werden aufgefordert, sich diese Filme anzusehen und herauszufinden, wie die Genialität, der Individualismus, der Erfolg oder das Scheitern und das Ausgegrenztsein der Mathematiker*innen dargestellt werden und inwiefern diese Darstellungen eben nicht mit der Realität übereinstimmen bzw. welche Einsichten über amerikanische Mythen, Ideologien (bzw. ideologische Konstrukte) und Wertvorstellungen wir daraus gewinnen können. Eine Zusammenstellung von Filmszenen, die im Unterricht herangezogen werden können, findet man unter: https://people.math.harvard.edu/~knill/math-movies/

Für die Arbeit mit audiovisuellen Texten eignet sich ein mehrstufiger Ansatz, da dieser ein differenzierteres Verständnis und eine kritischere Wahrnehmung eines Films oder einer bestimmten Szene ermöglicht, ohne das Filmerlebnis der Schüler*innen zu zerstören (Henseler et al., 2011; Thaler, 2012). Die Möglichkeit, dieselbe Szene mehr als einmal zu sehen, kann den Schüler*innen helfen, allmählich von einem eher allgemeinen zu einem detaillierteren audiovisuellen Verständnis zu gelangen, indem sie sich auf verschiedene Aspekte einer Szene konzentrieren und das Zusammenwirken von Bild, Ton, Kameraperspektive und Verwendung von Farben usw. untersuchen und erkennen, dass die Sinngebungsprozesse nicht zufällig sind (Viebrock, 2016).

Das Anschauen der Filme, Filmszenen oder Folgen einzelner Serien kann durch Fragen strukturiert werden, die die Schüler*innen beantworten müssen:
- Hat die Figur so geredet und ausgesehen, wie du sie dir vorgestellt hast?
- Wie werden Mathematik und Mathematiker*innen dargestellt und charakterisiert?
- Treiben die weiblichen Charaktere die Handlung und/oder die Mathematik (mathematische Entwicklungen und Erkenntnisse) voran? Oder unterstützen sie die männlichen Charaktere?
- Wie realistisch und glaubwürdig sind die Darstellungen der Frauen? Was kann man daraus über Vorurteile gegenüber Frauen und/oder gegenüber Mathematikerinnen schließen?

Angesichts positiv besetzter mathematischer Frauenfiguren in Filmen wie *Hidden Figures* mag man einen emanzipatorischen Trend konstituieren. Hier sei allerdings auf die bereits erwähnte Studie von Fellus et al. (2022) verwiesen, die zeigt, dass die Rezeption von Geschichten, die geschlechts- und gruppenspezifische Stereotypen ablehnen, andere Vorurteile dennoch verstärken kann. Wichtig ist daher vor allem, dass jegliche Formen von Stereotypisierungen explizit von den Lehrpersonen angesprochen und problematisiert werden, damit sich die Vorurteile nicht noch verstärken (Minas, 2018).

Darüber hinaus können die Schüler*innen dazu aufgefordert werden, sich über die mathematischen Ansätze und Methoden zu informieren, die von den Mathematiker*innen in diesen Filmen verwendet werden, um Probleme zu lösen. In Anlehnung an *Inception* könnten die Schüler*innen beispielsweise aufgefordert werden, sich mit den mathematischen Grundlagen der Topologie und Geometrie vertraut zu machen, die zur Manipulation der Traumwelt verwendet werden, während durch *The Imitation Game* das Interesse der Schüler*innen geweckt wird, die mathematischen Grundlagen der Kryptografie und Informatik besser verstehen zu wollen.

Die Fernsehserie *Big Bang Theory* (2007-2019) zeichnet sich nicht nur durch ihren speziellen Humor und ihre zahlreichen popkulturellen Referenzen aus, sondern ist die wohl bisher einzige Sitcom, die dazu beigetragen hat, das Interesse an Physik, Mathematik und Informatik zu fördern. *The Big Bang Theory* bedient dabei alle Klischees: Sheldon Cooper (gespielt von Jim Parsons) ist ein überhebliches Genie, der sich nicht zurückhält seine Meinung zu äußern. Er hat Schwierigkeiten sich auf romantische Beziehungen einzulassen und seine autistische Veranlagung zwingt ihn zur Einhaltung bestimmter Abläufe.

Es ist wichtig, im Unterricht zu diskutieren, wie im Genre der Sitcom (engl. situational comedy) stereotype Eigenschaften übertrieben werden, um den Charakteren Komik zu verleihen. Dass die Reduzierung auf bestimmte Merkmale zu einer Eindimensionalität führt, sollte ebenso problematisiert werden wie die Frage, ob diese Darstellung im Rahmen des komödiantischen Formats noch angemessen oder bereits problematisch ist. Dabei können Fragen wie die folgenden hilfreich sein: Was erfahren wir über die Figur des Sheldon Cooper, (a) durch seine eigenen Handlungen und Aussagen und (b) durch die Handlungen und Aussagen anderer? Welche Dimensionen seines Charakters bleiben unbeleuchtet? Nach einer solchen kritischen, diskursiv geführten Intervention ist es auch möglich, konkrete mathematische Beispiele aus der Serie zu übernehmen:
- **Rubiks Würfel:** In der Folge „Die Roboter-Manipulation" bzw. „The Robotic Manipulation" (4. Staffel, Folge 1) erklärt Sheldon, wie man Rubiks Würfel lösen kann und fordert seinen Kumpel Howard zu einem Duell heraus. Howard baut daraufhin einen Roboter bzw. programmiert ein Programm, das den Würfel für ihn löst. Diese Folge könnte für

- **Schrödingergleichung:** In der Folge „Die Parkplatz-Eskalation" bzw. „The Parking Spot Escalation" (Staffel 6, Folge 9) kommt es an besagtem Parkplatz zum Wettstreit zwischen Sheldon und Barry Kripke. Sheldon, als der ausgewiesenere Experte in Quantenmechanik, kann die Frage zu Schrödingergleichung schneller beantworten, erkennt jedoch, dass der Wettstreit und die Auseinandersetzung mit Barry Kripke sinnlos war. Diese Episode und Schrödingers Gleichung könnte als Einführung in die grundlegenden Konzepte der Quantenmechanik und Differentialgleichungen dienen.
- **Zahlentheorie:** Sheldons Leidenschaft für Zahlentheorie zeigt sich in mehreren Folgen, in denen Primzahlen eine Rolle spielen. In „Die Leuchtfisch-Idee" bzw. „The Luminous Fish Effect" (Staffel 1, Folge 4) arbeiten Sheldon und Howard an einem Forschungsprojekt, bei dem sie Kenntnisse der Primzahltheorie anwenden und eine Kamera entwickeln, die Leuchtfische im Dunkeln aufspüren kann. In Staffel 10, Folge 13 „Die Neuvermessung der Liebe" bzw. „The Romance Recalibration" nutzt Sheldon die Primzahltheorie, um ein Problem zu lösen, das er mit einem Valentinstagsgeschenk für seine Freundin Amy hat. Die Primzahltheorie, die in mindestens drei weiteren Folgen eine Rolle spielt, könnte als Einführung in die Konzepte der Arithmetik und Algebra dienen (siehe Staffel 3, Folge 10, „Die andere Seite der Krawatte"; Staffel 6, Folge 20, „Das Gürteltierproblem"; Staffel 7, Folge 5, „Die Krokodil-Summenformel").

In der Fachzeitschrift des National Council of Teachers of Mathematics findet man einige von den *Harry Potter*- und den *Hunger Games*-Filmen inspirierte mathematische Aktivitäten (Bush & Karp, 2012; Howe, 2002; Russo & Russo, 2017). Die folgenden Publikationen unterstützen Lehrkräfte bei der Auswahl und Einsatz populärer Medien:
- *Math, Culture, and Popular Media: Activities to Engage Middle School Students Through Film, Literature, and the Internet* (Chappell & Thompson, 2009).
- *Teaching Mathematics Using Popular Culture: Strategies for Common Core Instruction from Film and Television* (Reiser, 2015).

Die kommerziell von Texas Instruments unterstützte Webseite www.stemhollywood.com bietet ebenfalls eine Reihe von Aktivitäten für den Mathematikunterricht.

## Schlussbemerkung

Mathematik als Wissenschaft und die Mathematiker*innen leiden seit Jahrzehnten unter problematischen Wahrnehmungen in der Öffentlichkeit. Dabei steht außer Frage, dass Frauen mindestens so gute Mathematikerinnen sein können wie Männer, was nicht zuletzt durch die lange Liste bedeutender Mathematikerinnen belegt wird. Auch wenn es Beispiele für Kindheitstalente gibt, die sich zu guten Mathematiker*innen entwickelt haben, ist die Liste derjenigen die später „verschwinden" ebenso lang oder länger. Auf der anderen Seite hatten mehrere bedeutende Persönlichkeiten Lernschwierigkeiten als Kinder, wie zum Beispiel David Hilbert, der vielleicht einflussreichste Mathematiker des letzten Jahrhunderts. Gleichzeitig gibt es Beispiele für talentierte Mathematiker*innen, die im hohen Alter mit geistigen Defiziten zu kämpfen hatten. Die überwiegende Mehrheit der herausragenden Talente ist hingegen von ganz normalen Menschen mit Familie und Hobbies nicht zu unterscheiden. Mathematik als Wissenschaft wird in der Öffentlichkeit oft mit Rechnen gleichgestellt, während Theoriebilden, Problemlösen und Algorithmen entwickeln fast nie wahrgenommen wird. Es ist daher wichtig festzustellen, dass solche Beispiele allein nicht ausreichen, um Vorurteile abzubauen, sondern dass es die oben beschriebene Analyse stereotyper Darstellungen braucht, um eine tendenziöse Darstellung wahrnehmen zu können.

Es ist beunruhigend wie hartnäckig sich problematische Vorstellungen in unserer Gesellschaft halten und durch Mediendarstellungen verstärken können. Dass dieses falsche Bild die Ablehnung des Schulfachs Mathematik durch Schüler*innen begünstigt, wie die zitierte didaktische Forschung zeigt, ist umso alarmierender. Der kritische Blick auf Hollywoodproduktionen und Serienhits im Mathematik- und Englischunterricht erscheint uns als *eine* Möglichkeit diese Entwicklung sichtbar zu machen und ihr entgegenzutreten.

Unser Anliegen ist es, durch die Analyse von Blockbustern und populären Serien den Schüler*innen zu vermitteln wie sie verzerrte Darstellungen erkennen, Gründe benennen und – im Idealfall – Gegenentwürfe erstellen können. Dabei sei angemerkt, dass die Auswahl der zuvor diskutierten Filme und Serien rein exemplarischer Natur ist und sowohl unseren Vorlieben als auch Abneigungen geschuldet ist. Unsere Liste kann und sollte durch zahlreiche andere und aktuellere Produktionen erweitert werden. Da unser Beitrag dazu anregen möchte, den Mathematik- und Englischunterricht zu verzahnen, wäre es wünschenswert, auf Grundlage unserer Überlegungen eine Unterrichtssequenz zu entwickeln, diese zu erproben und empirisch mit einer Interventionsstudie und Interviews zu begleiten.

**Literatur**